\documentclass{article}

\usepackage{arxiv}

\usepackage[utf8]{inputenc} 
\usepackage[T1]{fontenc}    
\usepackage{hyperref}       
\usepackage{url}            
\usepackage{booktabs}       
\usepackage{amsfonts}       
\usepackage{nicefrac}       
\usepackage{microtype}      
\usepackage{graphicx}
\usepackage{natbib}
\usepackage{doi}
\usepackage{physics,amsmath}
\usepackage{authblk}

\title{Conservation and stability in a discontinuous Galerkin method for the vector invariant spherical shallow water equations}

\author[1]{Kieran Ricardo}
\author[2]{David Lee}
\author[1]{Kenneth Duru}
\affil[1]{Mathematical Sciences Institute, Australian National University, Canberra, Australia}
\affil[2]{Bureau of Meteorology, Melbourne, Australia}

\date{} 					

\begin{document}
\maketitle

\begin{abstract}
We develop a novel and efficient discontinuous Galerkin spectral element method (DG-SEM) for the spherical rotating shallow water equations in vector invariant form. We prove that the DG-SEM is energy stable, and discretely conserves mass,  vorticity, and linear geostrophic balance on general curvlinear meshes. These theoretical results are possible due to  our novel entropy stable numerical DG fluxes for the shallow water equations in vector invariant form. We experimentally verify these results on a cubed sphere mesh. Additionally, we show that our method is robust, that is can be run stably without any dissipation. The entropy stable fluxes are sufficient to control the grid scale noise generated by geostrophic turbulence without the need for artificial stabilisation. 
\end{abstract}

\section{Introduction}

Variational methods such as mixed finite element methods, and continuous Galerkin and discontinuous Galerkin spectral element methods (CG-SEM and DG-SEM respectively), have become an increasingly appealing choice for atmospheric modelling, owing to their scalable performance on modern supercomputers and ability to discretely preserve conservation properties \citep{mcrae2014energy, lee2018discrete, taylor2010compatible, gassner2016well, waruszewski2022entropy}. In this paper we present an entropy stable DG-SEM method for the rotating shallow water equations (RSWE) in vector invariant form that combines many of the strengths of these three approaches.

Ideally numerical methods should mimic the properties of the continuous system they approximate, but continuous systems have infinite invariants of which only a finite number can be preserved by any discretisation \citep{thuburn2008some}. This necessitates choosing a subset of the invariants, based on the particular modelling application, for the discrete model to conserve. For the RSWE in the context of atmospheric modelling we target the following: 
\begin{itemize}
	\item Local conservation of mass. Arguably the most important of these properties, it helps to limit the bias of long climate simulations \citep{thuburn2008some}. 
	\item  Local conservation of {absolute} vorticity. Mid-latitude weather is dominated by vorticity dynamics and its conservation assists in accurately modelling these dynamics \citep{lee2018discrete}. Additionally this prevents the gravitational potential, by far the largest component of atmospheric energy, from spuriously forcing the {absolute} vorticity which may destroy the meteorological signal \citep{staniforth2012horizontal}. 
 \item Local energy conservation and stability. The vast majority of energy dissipation within the atmosphere is due to three-dimensional effects not included in the RSWE. However there is evidence of a small amount of energy dissipation from two-dimensional layer wise effects on the order or $0.1Wm^{-2}$  \citep{thuburn2008some}. Additionally the RSWE conserves energy and so desirable discretisations of the RSWE should have minimal or no energy dissipation. Energy is also the entropy of the RSWE, and we follow the approach of using entropy stable fluxes to ensure non-linear stability \citep{gassner2016well,waruszewski2022entropy,carpenter2014entropy}. 
 \item Exact steady discrete geostrophic balance. The atmosphere is often characterised by slow moving variations around a linear state, consisting of a collection of slowly evolving geostrophic modes. To accurately represent this process a discretisation of the RSWE applied to the linear RSWE should contain exactly steady discrete geostrophic modes. \citep{cotter2012mixed,lee2018discrete}.
\end{itemize}
 
The {mixed finite element} approach is able to locally conserve mass, has steady discrete geostrophic modes, and globally conserves {absolute} vorticity and energy \citep{cotter2012mixed, mcrae2014energy, lee2018discrete}. {Mixed finite element} methods can also conserve potential enstrophy on affine geometries subject to exact integration of the chain rule, but this has not yet been extended to spherical geometry \citep{lee2018discrete, mcrae2014energy}.  Impressively, these schemes can be run stably without any dissipation \citep{lee2022comparison}. To simultaneously achieve all these properties {mixed finite element} methods use a non-diagonal mass matrix making them more computationally expensive than alternative spectral element methods \citep{lee2018discrete}. 

In parallel to this, a highly scalable collocated CG-SEM \citep{taylor2010compatible} was developed for the RSWE. CG-SEM applied to the vector invariant form of the RSWE has been shown to locally conserve mass, {absolute} vorticity, and energy. Due to dispersion errors CG-SEM requires filtering or artificial viscosity which may introduce undesirable numerical artefacts, excessively dissipate energy, and adversely effect the accuracy of the solution {\citep{melvin2012dispersion, dennis2012cam, ullrich2018impact}}.  

{The DG method has become a popular choice for modelling shallow water equations with many variants appearing in the literature \citep{giraldo2002nodal, nair2005discontinuous, gassner2016well, wintermeyer2017entropy, dumbser2013staggered}. A DG method for the RSWE was introduced in \citep{giraldo2002nodal}, this approach used an icosahedral grid and Cartesian coordinates. Following this {\citep{nair2005discontinuous}} developed the first DG method for the RSWE using a cubed sphere grid. In an alternate approach \citep{dumbser2013staggered} presents a DG method with staggered cells for the non-rotating shallow water equations. In \citep{gassner2016well} and \citep{wintermeyer2017entropy} the authors develop the first energy conserving DG method for the non-rotating shallow water equations using a split-form DG-SEM. In contrast to the split-form approach of \citep{gassner2016well} and \citep{wintermeyer2017entropy}, \citep{gaburro2023high} presents an entropy conserving DG method for general systems of conservation laws by adding correction terms to cancel spurious numerical entropy generation, this builds upon the work of \citep{abgrall2018general, abgrall2022reinterpretation} which develops a general framework for developing entropy conserving methods using entropy correcting terms. To the best of the authors knowledge, all DG methods for the  RSWE  presented in the literature discretise the conservative form of the equations and therefore do not conserve absolute vorticity or preserve linear geostrophic balance.} This is a direct result of solving the RSWE in conservative form as opposed to the vector invariant form. Similar to CG-SEM, DG-SEM possess dispersion errors and therefore requires that high frequency waves be damped. However, DG-SEM can achieve this damping through entropy stable numerical fluxes which do not affect the high order accuracy of the scheme.

{
In an alternate branch of research, \citep{castro2008finite} and \citep{castro2017well} develop well-balanced finite volume methods for the spherical shallow water equations without rotation. By design these methods accurately maintain steady state solutions of the shallow water equations but are more applicable to tsunami simulations where the Coriolis force is negligible.}

In this work we present a DG-SEM for the RSWE in vector invariant form that is arbitrarily high order accurate. We prove that the DG-SEM locally conserves mass and {absolute} vorticity, locally semi-discretely conserves energy, and preserves linear geostrophic balance. We then show how to consistently modify this scheme to dissipate energy/entropy through the use of dissipative fluxes. Numerical experiments are presented to verify the theoretical results and show that geostrophic turbulence can be well simulated with our entropy stable fluxes with no additional stabilisation strategy.

The rest of the article is as follows. Section 2 presents the continuous equations and derives the conservation and stability properties of the continuous model that should be  preserved by our numerical method. Section 3 introduces our spatial discretisation and derives several identities for the function spaces and their discrete operators. In section 4  we present theoretical results for discrete conservation and stability properties of the method. Section 5 presents the results of numerical experiments that verify our theoretical analysis. In section 6, we draw conclusions and suggests directions for future work.

\section{Continuous shallow water equations}

In this section we introduce the RSWE on the surface of a spherical geometry and derive the continuous properties that our numerical scheme should discretely mimic.

\subsection{Vector invariant equations}

We use the vector invariant form of the RSWE on a 2D manifold $\Omega$ embedded in $R^3$ with periodic boundary conditions. We choose periodic boundary conditions as they simplify the analysis and the domain of most practical interest is the sphere. The RSWE are 
\begin{equation}
	\vb{u}_t + \omega\vb{k}\cross\vb{u} + \nabla G = 0\text{,}
	\label{eq:momentum}
\end{equation}
\begin{equation}
	D_t + \nabla \cdot \vb{F} = 0\text{,}
	\label{eq:mass}
\end{equation}
with
\begin{equation}
	\omega = \vb{k} \cdot \nabla \cross \vb{u} + f\text{,} \quad \vb{F} = D\vb{u}\text{,} \quad G = \tfrac{1}{2}\vb{u}\cdot\vb{u} + gD\text{,}
\end{equation}
where $t\ge 0$ denotes the time variable, the unknowns are the flow velocity vector $\vb{u}=[u, v, w]^T$  and the fluid depth $D$. {Note that here and in the rest of the paper all differential operators are taken to be constrained to the manifold.} Furthermore, $\omega$ is the {absolute} vorticity, $f$ is the Coriolis frequency, $\vb{F}$ is the mass flux, $G$ is the potential, and $\vb{k}$ is normal of the manifold. At $t=0$, we augment the RSWE \eqref{eq:momentum}--\eqref{eq:mass} with the initial conditions
\begin{align}
    \vb{u} = \vb{u}_0(\mathbf{x}), \quad {D} = {D}_0(\mathbf{x}), \quad \mathbf{x} \in \Omega,
\end{align}
such that the flow is constrained to the manifold by ensuring that $\vb{u}_0 \cdot \vb{k} = 0$. {The constraint of the differential operators and the initial constraint $\vb{u}_0 \cdot \vb{k} = 0$ is sufficient to ensure that the $\vb{u}(\vb{x}, t) \cdot \vb{k} = 0$ holds $\forall \vb{x}\in \Omega,\; t \geq 0$.}

\subsection{Conservation of mass, {absolute} vorticity, and energy}\label{sec:Conserve_mass_vorticity_energy}
We will show that the total mass, {absolute} vorticity and energy are conserved. Integrating the continuity equation \eqref{eq:mass} over a periodic domain yields conservation of total mass, that is
\begin{align}\label{eq:mass_conservation}
     \frac{d}{dt}\int_{\Omega} D d\Omega = 0.
\end{align}
We can derive a continuity equation for {absolute} vorticity $\omega$ by applying $\vb{k} \cdot \nabla \cross$ to the momentum equation \eqref{eq:momentum} giving
\begin{equation}\label{eq:vorticity_continuity}
	\omega_t + \nabla \cdot (\omega\vb{u}) = 0\text{.}
\end{equation}
Similarly, integrating equation \eqref{eq:vorticity_continuity} over a periodic domain yields conservation of total {absolute} vorticity,
\begin{align}\label{eq:vorticity_conservation}
      \frac{d}{dt}\int_{\Omega} \omega d\Omega = 0.
\end{align}
{Additionally, all quantities of the form $Db\left(\frac{\omega}{ D}\right)$, for arbitrary function $b$ of $\frac{\omega}{ D}$, are also locally conserved by the RSWE \citep{vallis2017atmospheric}. The potential vorticity $\frac{\omega}{D}$ is transported by the flow, therefore $b\left(\frac{\omega}{ D}\right)$ is also transported by the flow, and hence $Db\left(\frac{\omega}{ D}\right)$ is locally conserved. The potential enstrophy is one such quantity, and while we do not target its conservation here several other numerical schemes have \citep{mcrae2014energy, lee2018discrete}}.

Let the elemental energy $E$ be denoted by
\begin{equation}
   E=\tfrac{1}{2}D\vb{u}\cdot\vb{u} + \tfrac{1}{2}gD^2.
\end{equation}
The elemental energy $E$ also satisfies a continuity equation. To show this we take the time derivative of the elemental energy $E$ giving
\begin{equation}
     E_t = \vb{F}\cdot \vb{u}_t + (\tfrac{1}{2}\vb{u}\cdot\vb{u} + gD)D_t   = -\vb{F}\cdot \omega \vb{k} \cross \vb{u} - \vb{F}\cdot \nabla G - G\nabla \cdot \vb{F}.
 \end{equation}
 We have the continuity equation for the energy $E$
\begin{equation}\label{eq:energy_continuity}
    E_t + \nabla \cdot \big(G\vb{F}\big) = 0,
\end{equation}
where we have used $\vb{F}\cdot \omega \vb{k} \cross \vb{u} = \omega D \vb{u}\cdot  \vb{k} \cross \vb{u} = 0$ and the product rule for derivatives.
As before, integrating equation \eqref{eq:energy_continuity} over a periodic domain yields conservation of total energy,
\begin{align}\label{eq:energy_conservation}
 \frac{d}{dt}\int_{\Omega} E d\Omega = 0.
\end{align}

\subsection{Linear geostrophic balance}\label{sec:linear_geo_balance}

{On domains with constant Coriolis parameter $f$}, in addition to a set of inertio-gravity wave solutions, the linear rotating shallow water equations also exhibit steady geostrophic modes \citep{cotter2012mixed, vallis2017atmospheric}. The RSWE \eqref{eq:momentum}--\eqref{eq:mass} linearised around a state of constant mean height $H$ and zero mean flow are
\begin{equation}
	\vb{u}_t + f{\vb{k} \cross} \vb{u} + g\nabla D = 0\text{,}
\end{equation} 
\begin{equation}
	D_t + H\nabla \cdot \vb{u} = 0\text{.}
\end{equation}
The steady solutions to these equations can be found by choosing an arbitrary stream function $\psi$ and setting
$D = -\tfrac{f}{g}\psi, \quad \vb{u} = \nabla {\cross \psi \vb{k}}$. 
Then we have
\begin{equation}
    D_t = -{H}\nabla \cdot \vb{u} = -{H\nabla \cdot \nabla \cross \psi \vb{k}} = 0,
\end{equation}
\begin{equation}
    \vb{u}_t = -(f {\vb{k} \cross }\vb{u}  + g\nabla D) = -(f\nabla \psi - f\nabla \psi{) = 0}.
\end{equation}

\subsection{Entropy}

 An entropy function is any convex combination of prognostic variables that also satisfies a conservation law \citep{leveque1992numerical}. The importance of entropy functions is that they are conserved in smooth regions and are dissipated across discontinuities in physically relevant weak solutions. DG methods are continuous within each element but can be discontinuous across element boundaries, this motivates that the volume terms of DG methods should conserve entropy while numerical fluxes should dissipate entropy. We call this local entropy stability and in practice many numerical methods have found that this is sufficient for numerical stability \citep{giraldo2020introduction}.
 
 It is well known that energy is an entropy function of the RSWE in conservative form \citep{gassner2016well}, here we show that energy is also an entropy for the RSWE in vector invariant form for subcritical flow $|\vb{u}| < \sqrt{gD}$. Atmospheric flows are subcritical and so this restriction does not cause any difficulty in practice, for the target applications. The analysis in section 2.3 shows that energy is conserved, all that remains is to show that energy is also a convex combination of $u$, $v$, $w$, and $D$. {In the following analysis we non-dimensionalise all quantities by their mean values to avoid issues with conflicting units.} The Hessian of the energy with respect to the prognostic variables is
\begin{equation}
    H_E = \begin{pmatrix}
    D &  0  & 0 & u \\
    0 &  D  & 0 & v \\
     0 & 0  & D & w \\
    u &  v  & w & g \\
    \end{pmatrix}\text{,}
\end{equation}
and its characteristic equation is
\begin{equation}
    (D-\lambda)^2\big(\lambda^2 - (D+g)\lambda + Dg - |\vb{u}|^2\big)\text{.}
\end{equation}
As $D,g>0$ the eigenvalues corresponding to the left bracket are always positive, and as long as $Dg - |\vb{u}|^2 > 0$, which implies $|\vb{u}| < \sqrt{gD}$, the eigenvalues corresponding to the quadratic in the second bracket are also positive. Therefore for subcritical flow the energy is an entropy function of the vector invariant RSWE.

An effective numerical method should as far as possible conserve mass, {absolute} vorticity, energy, entropy and linear geostrophic balance.
\section{Spatial discretisation}

In this section we describe the spatial discretisation of our method and prove several discrete identities that are necessary for numerical stability and discrete conservation properties. Our method decomposes the domain into quadrilateral elements $\Omega^q$ and approximates the solution in each quadrilateral by Lagrange polynomials. The solutions are connected across the elements through energy-stable numerical fluxes.

\subsection{Function spaces}

For each element $\Omega^q$ we define an invertible map $r(\vb{x}; q)$ to the reference square $[-1, 1]^2$, and approximate solutions by polynomials of order $m$ in this reference square. We use a computationally efficient tensor product Lagrange polynomial basis with interpolation points collocated with Gauss-Lobatto-Legendre (GLL) quadrature points. Using this basis the polynomials of order $m$ on the reference square can be defined as
\begin{equation}
    P_m := \text{span}_{i, j=1}^{m+1} l_i(\xi) l_j(\eta)\text{,}
\end{equation}
where $(\xi, \eta) \in [-1, 1]^2$ and $l_i$ is the 1D Lagrange polynomial which interpolates the $i^{th}$ GLL node. We define a discontinuous scalar space $S$ which contains the depth $D^h$ as
\begin{equation}
    S = \{\phi \in L^2(\Omega) : r(\phi; q) \in P_m, \forall q \} \text{.}
\end{equation}
Note that within each element functions $\phi \in S$ can be expressed as
\begin{equation}
    \phi(\vb{x}) = \sum_{i,j=1}^{m+1}{\phi_{ijq}l_i(\xi)l_j(\eta)}, \forall \vb{x} \in \Omega^q\text{,}
\end{equation}
where $\xi, \eta = r(\vb{x}; q)$, $\phi_{ijq} = \phi(r^{-1}(\xi_i, \eta_j;q))$.

We define a discontinuous vector space which contains the velocity $\vb{u}$. Let $\vb{v}_1(x,y,z)$ and $\vb{v}_2(x,y,z)$ be any independent set of vectors which span the tangent space of the manifold at each point $(x,y,z)$, then the discontinuous vector space $V$ containing the velocity $\vb{u}$ can be defined as
\begin{equation}
    V = \{\vb{w} : \vb{w} \cdot \vb{v}_i \in S, i=1,2 \}\text{.}
\end{equation}
We note that $V$ is independent of the particular choice of $\vb{v}_i$ however two convenient choices are the contravariant and covariant {basis} vectors. 

\subsection{Covariant and contravariant {base} vectors}

Here we provide a brief summary of the differential geometry tools we use in our discrete method. This summary closely follows \citep{taylor2010compatible} and these results can also be found in standard textbooks \citep{giraldo2020introduction, heinbockel2001introduction}. The covariant {base} vectors $\vb{g}_1$ and $\vb{g}_2$ are defined as
\begin{equation}
    \vb{g}_1 = \dfrac{\partial \vb{x}}{\partial \xi}, \quad \vb{g}_2 = \dfrac{\partial \vb{x}}{\partial \eta}\text{,}
\end{equation}
with these the determinant of the Jacobian of $r(\vb{x}; n)$ can be expressed as
\begin{equation}
    J = |\vb{g}_1 \times \vb{g}_2|\text{.}
\end{equation}
Further, we assume that the map $r(\vb{x}; q)$ is invertible and the Jacobian is strictly positive $J>0$.

The contravariant {base} vectors are {defined as}
\begin{equation}
    \vb{g}^1 = \nabla \xi, \quad {\vb{g}^2} = \nabla \eta \text{}
\end{equation}
{to fulfill the property}
\begin{equation}
    \vb{g}^\alpha \cdot \vb{g}_\beta = \delta{^\alpha_\beta}\text{,}
\end{equation}
enabling vectors to {be} expressed as
\begin{equation}
    \vb{w} = \sum_{\alpha=1,2}{w^\alpha\vb{g}_\alpha} = \sum_{\alpha=1,2}{w_\alpha\vb{g}^\alpha}\text{,}
\end{equation}
where $w^\alpha = \vb{w}\cdot \vb{g}^\alpha$, $w_\beta = \vb{w}\cdot \vb{g}_\beta$.

For flows such that $\vb{w}\cdot\vb{k}=0$ this enables the divergence, gradient, and curl to be expressed as
\begin{equation}
    \nabla \cdot \vb{w} = \dfrac{1}{J}\bigg(\dfrac{\partial Jw^1}{\partial \xi} + \dfrac{\partial Jw^2}{\partial \eta} \bigg)\text{,}
    \label{eq:div}
\end{equation}
\begin{equation}
    \nabla \phi = \dfrac{\partial \phi}{\partial \xi} \vb{g}^1 + \dfrac{\partial \phi}{\partial \eta} \vb{g}^2\text{,} 
    \label{eq:grad}
\end{equation}
\begin{equation}
    \nabla \cross \vb{w} = \dfrac{1}{J}\bigg(\dfrac{\partial w_2}{\partial \eta} - \dfrac{\partial w_1}{\partial \xi} \bigg)\vb{k}\text{,} 
    \label{eq:curl}
\end{equation}
\begin{equation}
    \nabla \cross \phi \vb{k} = \dfrac{1}{J}\dfrac{\partial \phi}{\partial \eta}\vb{g}_1 - \dfrac{1}{J}\dfrac{\partial \phi}{\partial \xi} \vb{g}_2 \text{,} 
    \label{eq:curl2}
\end{equation}
where we have only used $\nabla \cross$ on quantities either parallel or orthogonal to $\vb{k}$ to simplify the exposition.

\subsection{Discrete inner product and boundary integrals}

We approximate integrals over the elements by using GLL quadrature. This defines a discrete element inner product
\begin{equation}
    \langle f, g \rangle_{\Omega^q} := \sum_{i, j=1}^{m+1}{w_i w_j J_{ij} f_{ij}g_{ij}}\text{,}
\end{equation}
which approximates the integral
\begin{equation}
    \int_{\Omega^q}{fg d\Omega^q} \approx \langle f, g \rangle_{\Omega^q}\text{.}
\end{equation}
The discrete global inner product is defined simply as $\langle f, g \rangle_{\Omega} = \sum_q{\langle f, g \rangle_{\Omega^q}}$. Similarly for vectors 
\begin{equation}
    \langle \vb{f}, \vb{g} \rangle_{\Omega^q} := \sum_{i, j=1}^{m+1}{w_i w_j J_{ij} \vb{f}_{ij} \cdot \vb{g}_{ij}}\text{.}
\end{equation}
We also approximate element boundary integrals using GLL quadrature
\begin{equation}
    \int_{\partial \Omega^q}{\vb{f}\cdot\vb{n}dl} \approx \langle \vb{f}, \vb{n}\rangle_{\partial\Omega^q}\text{,}
\end{equation}
where 
\begin{equation}
\begin{split}
    \langle \vb{f}, \vb{n}\rangle_{\Omega^q} := \sum_{j}{w_j\big(|g_1|_{1j} \vb{f}_{1j}\cdot \vb{n}_{1j} + |g_1|_{(m+1)j} \vb{f}_{(m+1)j}\cdot \vb{n}_{(m+1)j}} \\ {+ |g_2|_{j1}\vb{f}_{j1}\cdot \vb{n}_{j1} + |g_2|_{j(m+1)}\vb{f}_{j(m+1)}\cdot \vb{n}_{j(m+1)}\big)}\text{,}
\end{split}
\end{equation}
and $\vb{n}$ is the outward facing unit normal of the element boundary $\partial \Omega^q$ given by
\begin{equation}
\vb{n} = 
    \begin{cases}
        \xi \dfrac{\vb{g}^1}{|\vb{g}^1|}, \;\text{for}\; \xi = -1, 1 \\ 
        \eta \dfrac{\vb{g}^2}{|\vb{g}^2|}, \;\text{for}\; \eta = -1, 1\text{.}
\end{cases}
\end{equation}
The unit tangent vector $\vb{t}$ is similarly defined as 
\begin{equation}
\vb{t} = 
    \begin{cases}
        \xi \dfrac{\vb{g}_2}{|\vb{g}_2|}, \;\text{for}\; \xi = -1, 1 \\ 
        -\eta \dfrac{\vb{g}_1}{|\vb{g}_1|}, \;\text{for}\; \eta = -1, 1
\end{cases}\text{.}
\end{equation}

With the discrete inner products defined we now construct discrete $L^2$ projection operators $\Pi_S(\cdot )$ and $\Pi_V(\cdot )$ which preserve discrete $L^2$ inner products in $S$ and $V$ respectively. For example $\Pi_S(\cdot )$ is defined as the unique projection that satisfies
\begin{equation}
	\langle \gamma, b \rangle_\Omega = \langle \gamma, \Pi_S(b) \rangle_\Omega\text{,} \; \forall \gamma \in S\text{,}
\end{equation} 
for all functions $b$. $\Pi(\cdot )_V$ is defined similarly. In practice the discontinuous projection operators $\Pi_S$ and $\Pi_V$ simply interpolate functions within each element
\begin{equation}
    \Pi_S(f)(\vb{x}) = \sum_{ij}{f_{ijq}l_i(\xi)l_j(\eta)}, \forall \vb{x} \in \Omega^q\text{,}
\end{equation}
\begin{equation}
    \Pi_V(\vb{f})(\vb{x}) = \sum_{ij}{\vb{f}_{ijq}l_i(\xi)l_j(\eta)}, \forall \vb{x} \in \Omega^q\text{,}
\end{equation}
where $(\xi, \eta) = r(\vb{x}; q)$.

\subsection{Discrete operators}

Following \citep{taylor2010compatible} we construct discrete operators by discretising equations \eqref{eq:div}--\eqref{eq:curl2}
\begin{equation}
    \nabla_d \cdot \vb{w} = \Pi_S\Bigg(\dfrac{1}{J}\bigg(\dfrac{\partial}{\partial \xi} \Pi_S(Jw^1) + \dfrac{\partial}{\partial \eta} \Pi_S(Jw^2)\bigg)\Bigg)\text{.}
    \label{eq:div-disc}
\end{equation}
\begin{equation}
    \nabla_d \phi = \dfrac{\partial \phi}{\partial \xi}  \vb{g}^1 + \dfrac{\partial \phi}{\partial \eta} \vb{g}^2\text{,} 
    \label{eq:grad-disc}
\end{equation}
\begin{equation}
    \nabla_d \cross \vb{w} = \Pi_S\Bigg(\dfrac{1}{J}\bigg(\dfrac{\partial w_2}{\partial \xi} -\dfrac{\partial w_1}{\partial \eta} \bigg)\Bigg)\vb{k} \text{,} 
    \label{eq:curl-disc}
\end{equation}
\begin{equation}
    \nabla_d \cross \phi \vb{k} = \Pi_V\Bigg(\dfrac{1}{J}\bigg(\dfrac{\partial \phi}{\partial\eta} \vb{g}_1-\dfrac{\partial \phi}{\partial\xi} \vb{g}_2\bigg)\Bigg)\text{,}
\end{equation}

\subsection{Properties of function spaces}

Here we present the properties of the function spaces and discrete operators necessary for the numerical stability and discrete conservation proofs performed later in this paper.

\subsubsection{Element local summation by parts}

The discrete gradient and divergence operators satisfy an element local summation by parts (SBP) property {\citep{taylor2010compatible}}
\begin{equation}
	\langle \nabla_d \phi, \vb{w} \rangle_{\Omega^q} + \langle \phi, \nabla_d \cdot \vb{w} \rangle_{\Omega^q} = \langle \phi \vb{w}, \vb{n}\rangle_{\partial \Omega^q}\text{,}\; \forall \phi \in S\text{,}\; \vb{w} \in V\text{.}
\end{equation}
This allows us to transform between the weak and strong form of our DG prognostic equations which we use in several proofs.

\subsubsection{Divergence curl cancellation}
The discrete curl operator satisfies a discrete $grad-curl$ cancellation
\begin{equation}
	\nabla_d \cdot \nabla_d \cross \phi \vb{k} = 0, \forall \phi \in S\text{,}
\end{equation}
which we use for {absolute} vorticity conservation. 

\subsubsection{Continuity of interpolants}

For tensor product bases with boundary interpolation points, the element local interpolation of continuous functions $\psi^h = \Pi_S(\psi)$ are continuous across element boundaries. To see this consider two elements $q_1$ and $q_2$ that share a boundary at $\xi=1$ and $\xi=-1$ respectively. The value of the interpolant $\psi^h$ in element $q_1$ at the shared boundary is
\begin{equation}
    \psi^h(1, \eta; q_1) = \sum_{i, j=1}^{m+1}{\psi_{ijq_1}l_i(1)l_j(\eta)} = \sum_{j=1}^{m+1}{\psi_{(m+1)jq_1}l_j(\eta)}\text{,}
\end{equation}
and in element $q_2$ at that same boundary
\begin{equation}
    \psi^h(-1, \eta; q_2) = \sum_{i, j=1}^{m+1}{\psi_{ijq_2}l_i(-1)l_j(\eta)} = \sum_{ j=1}^{m+1}{\psi_{1jq_2}l_j(\eta)}\text{.}
\end{equation}
The mappings $r(\cdot; q)$ agree along element boundaries, therefore $r^{-1}(1, \eta; q_1) = r^{-1}(-1, \eta; q_2)$ and
\begin{equation}
    \psi_{(m+1)jk_1} = \psi\big(r^{-1}(1, \eta_j; q_1)\big) = \psi\big(r^{-1}(-1, \eta_j; q_2)\big) = \psi_{1jk_2}\text{,}
\end{equation}
implying that $\psi^h$ is continuous across element boundaries.

\subsubsection{Curl gradient identity}

We also note that 
\begin{equation}
    \nabla_d \cross \phi\vb{k} = -{\Pi_V\big(}\vb{k} \cross \nabla_d \phi{\big)}, \; \forall \phi \in S\text{,}
    \label{eq:curl-grad-k}
\end{equation}
which we use in our discrete linear geostrophic balance proof.

\subsubsection{Element boundary summation by parts}\label{sec:element_boundary_sbp}

The final property we use is an element boundary SBP 
\begin{equation}
    \langle \phi \nabla_d \gamma + \gamma \nabla_d \phi, \vb{t} \rangle_{\partial \Omega^q} = 0\text{,}\; \forall \gamma, \phi \in S\text{.}
\end{equation}
This arises from a 1D SBP property for integrals along each side of the quadrilateral. Consider the $\xi=1$ boundary, $|\vb{g_2}|\nabla_d \phi \cdot \vb{t} = \dfrac{\partial \phi}{\partial \eta}$ and therefore by the exact integration of $2m-1$ polynomials with GLL quadrature
\begin{equation}
\begin{split}
    \langle \phi \nabla_d \gamma + \gamma \nabla_d \phi, \vb{t} \rangle_{\partial \Omega^q} = \sum_{\xi=-1, 1}\int_{\eta=-1}^{\eta=1}{\phi\dfrac{\partial \gamma}{\partial \eta} + \gamma\dfrac{\partial \phi}{\partial \eta} d\eta} \\ + \sum_{\eta=-1, 1}\int_{\xi=-1}^{\xi=1}{\phi\dfrac{\partial \gamma}{\partial \eta} + \gamma\dfrac{\partial \phi}{\partial \eta} d\xi} = 0\text{.}
\end{split}
\end{equation}

\subsection{The semi-discrete formulation}

First we introduce notation for jumps $[[\cdot]]$ and averages $\{\{\cdot \}\}$ across element boundaries as
\begin{equation}
    \{ \{ a \}\} = \tfrac{1}{2}\big(a^{+} + a^{-}\big)\text{,}
\end{equation}
\begin{equation}
    [[ a ]] = a^{+} - a^{-}\text{,}
\end{equation}
where $\cdot^{+}$ and $\cdot^{-}$ denote points on opposites sides of the boundary element boundary, and $a$ can be either a scalar or vector. We now define our spatial discretisation as
\begin{equation}
\begin{split}
	\langle \vb{w}^h, \vb{u}_t^h\rangle_{\Omega^q} + \langle \vb{w}^h, \omega^h \vb{k}\cross \vb{u}^{h}\rangle_{\Omega^q}+\langle \vb{w}^h, \nabla_d G^h\rangle_{\Omega^q}  \\
+ \langle \vb{w}^h, (\hat{G}-G)\vb{n}\rangle_{\partial\Omega^q} = 0\text{,}\;\forall \vb{w}^h\in V\text{,}
\end{split}
\label{eq:dg-velocity}
\end{equation}
\begin{equation}
	\langle \phi^h, D_t^h\rangle_{\Omega^q} + \langle \phi^h, \nabla_d \cdot \vb{F}^h\rangle_{\Omega^q} + \langle \phi^h, \big(\hat{\vb{F}}-\vb{F}^h\big)\cdot \vb{n}\rangle_{\partial\Omega^q} = 0\text{,}\;\forall \phi^h\in S\text{,}
 \label{eq:dg-depth}
\end{equation}
where the numerical fluxes are defined for the scalars $\alpha,\beta\ge 0$ as
\begin{equation}\label{eq:potential_flux}
    \hat{G} = \{\{ G \}\} + \alpha[[\vb{F}^h]]\cdot{\vb{n}^{+}}\text{.}
\end{equation}
\begin{equation}\label{eq:mass_flux}
    \hat{\vb{F}} \cdot \vb{n}^{+} = \{\{ \vb{F}^h \}\}\cdot \vb{n}^{+} + \beta[[G^h]]\text{,}
\end{equation}
and the discrete {absolute} vorticity is 
\begin{equation}
    \langle \phi^h, \omega^h\rangle_{\Omega^q} = \langle -\nabla_d \cross \phi^h\vb{k}, \vb{u}^h\rangle_{\Omega^q} + \langle \phi^h, f\rangle_{\Omega^q} + \langle \phi^h, \{\{\vb{u}^h\}\} \cdot \vb{t}\rangle_{\partial\Omega^q}\;\forall \phi^h\in S\text{,}
    \label{eq:dg-vort}
\end{equation}
where $\vb{t}$ is the unit tangent vector along element boundaries.

As will be shown below, our scheme is discretely entropy stable so long as $\alpha, \beta \geq 0$, and is locally semi-discretely energy/entropy conserving  for $\alpha = \beta = 0$ and locally discretely energy/entropy dissipating for $\alpha, \beta > 0$. We use a centred mass flux $\beta = 0$ to ensure that potential energy is not dissipated, and either use centred fluxes $\alpha=0$ or {  
 \begin{equation}
 \alpha = \frac{1}{2}\max\left(\frac{c^+}{D^+}, \frac{c^-}{D^-}\right)\text{,}
  \end{equation}
  where $c^\pm =|\vb{u}^{\pm}| + \sqrt{gD^{\pm}}$} is the fastest wave speed on either side of the boundary. This ensures that $\hat{G}$  has the correct units and reduces to a Rusanov flux \citep{rusanov1970difference} in the case of the linearised RSWE in section 2.5. We refer to the parameter choices $\alpha=\beta=0$ and $\alpha=\frac{1}{2}\max\left(\frac{c^+}{D^+}, \frac{c^-}{D^-}\right),\;\beta=0$ as the energy conserving and energy dissipating methods.

\section{Semi-discrete conservation and numerical stability}
In this section we prove some discrete conservation properties and establish the numerical stability of the numerical method \eqref{eq:dg-velocity}--\eqref{eq:dg-vort}. The main results of this paper are presented in this section, including conservation  of mass, {absolute} vorticity,  entropy and linear geostrophic balance at the discrete level, and numerical energy stability.

\subsection{Conservation of mass}

Element local conservation of total mass can be shown using the element local SBP property to transform the mass equation (\ref{eq:mass}) into weak form
\begin{equation}
	\langle \phi^h, D_t^h\rangle_{\Omega^q} - \langle \nabla_d \phi^h,  \vb{F}^h\rangle_{\Omega^q} + \langle\phi, \hat{\vb{F}}\cdot \vb{n}\rangle_{\partial\Omega^q} = 0\text{,}\;\forall \phi^h\in S\text{,}
	\label{eq:weak-mass}
\end{equation}
and then choosing the test function $\phi^h = 1$ to yield the  rate of change of the total mass of element $k$
\begin{equation}
	\langle 1, D_t^h\rangle_{\Omega^q} + \langle 1, \hat{\vb{F}}\cdot \vb{n}\rangle_{\partial\Omega^q} = 0\text{.}
	\label{eq:total-weak-mass}
\end{equation}
By construction the numerical fluxes are continuous across element boundaries and so the mass flowing through each point along the element boundary is equal and opposite to that of the adjacent element, therefore mass is locally conserved. Global conservation can be trivially obtained from local conservation in a periodic domain by noting that mass flux cancels for adjacent element and so $\sum_q{\langle 1, D_t^h\rangle_{\Omega^q}} = -\sum_q{\langle 1, \hat{\vb{F}}\cdot \vb{n}\rangle_{\partial\Omega^q}} = 0$.

\subsection{Conservation of {absolute} vorticity}

We show that the potential flux $G$  does not spuriously force the {absolute} vorticity, and that the {absolute} vorticity is also locally conserved. This is a direct consequence of our choice of the numerical flux $\hat{G}$ defined in \eqref{eq:potential_flux}, numerical approximation of the {absolute} vorticity $\omega^h$ giving in \eqref{eq:dg-vort}, and the element boundary SBP property discussed in subsection \ref{sec:element_boundary_sbp}.

The evolution of discrete {absolute} vorticity follows from taking the time derivative of $\omega^h$ in \eqref{eq:dg-vort}, giving
\begin{equation}
    \langle \phi^h, \omega^h_t\rangle_{\Omega^q} = -\langle \nabla_d \cross \phi^h\vb{k}, \vb{u}^h_t\rangle_{\Omega^q} + \langle \phi^h,  \{\{\vb{u}^h_t\}\} \cdot \vb{t}\rangle_{\partial\Omega^q}\text{,} \;\forall \phi^h\in S\text{.}
    \label{eq:dg-vort-evo}
\end{equation}
First we expand the volume term by using the weak form of the velocity equation \eqref{eq:dg-velocity} and the curl gradient identity \eqref{eq:curl-grad-k}, that is $\vb{n} \cdot \nabla_d \cross \phi^h \vb{k} = \nabla_d \phi^h \cdot \vb{t}$, and $\nabla_d \cdot \nabla_d \cross \phi^h\vb{k} = 0$ yielding
\begin{equation}
    \langle \nabla_d \cross \phi^h\vb{k}, \vb{u}^h_t\rangle_{\Omega^q} = \langle \nabla_d \phi^h, \omega^h \vb{u}^h\rangle_{\Omega^q} - \langle \hat{G}\nabla \phi^h,  \vb{t}\rangle_{\partial\Omega^q}\text{,} \;\forall \phi^h\in S\text{.}
\end{equation}
Next we evaluate the boundary term, as this is a GLL collocated method
\begin{equation}
    \langle l_i(\xi) l_j(\eta) \vb{t}_{ij}, \vb{u}^h_t \rangle_{\Omega^q} = w_i w_j (\vb{u}^h_{ij})_t \cdot \vb{t}_{ij}\text{,}
\end{equation}
substituting the test function $l_i(\xi)l_j(\eta) \vb{t}_{ij}$ into \eqref{eq:dg-velocity} yields
\begin{equation}
     \langle l_i(\xi) l_j(\eta) \vb{t}_{ij}, \vb{u}^h_t \rangle_{\Omega^q} = -w_i w_j \omega^h \vb{k} \cross \vb{u}^h\cdot \vb{t} - w_i w_j \nabla_d G\cdot \vb{t}\text{} 
\end{equation}
\begin{equation}
    = -w_i w_j \omega^h \vb{u}^h\cdot \vb{n} - w_i w_j\nabla_d G\cdot \vb{t}\text{,}
\end{equation}
and therefore
\begin{equation}
    \{\{\vb{u}_t\}\}\cdot \vb{t} = -\{\{\omega^h \vb{u}^h\}\} \cdot \vb{n} - \nabla_d\{\{ G\}\} \cdot \vb{t}\text{.}
\end{equation}
With these the {absolute} vorticity evolution becomes
\begin{equation}
\begin{split}
    \langle \phi^h, \omega^h_t\rangle_{\Omega^q} &= \langle \nabla_d \phi^h, \omega^h \vb{u}^h\rangle_{\Omega^q} - \langle  \hat{G}\nabla \phi^h + \phi^h \{\{\nabla_d G\}\},  \vb{t}\rangle_{\partial\Omega^q} \\
    &-\langle \phi^h \{\{\omega^h \vb{u}^h\}\},\vb{n}\rangle_{\partial\Omega^q} \\
    &=\langle \nabla_d \phi^h, \omega^h \vb{u}^h\rangle_{\Omega^q} - \langle  \{\{G\}\}\nabla \phi^h + \phi^h \{\{\nabla_d G\}\},  \vb{t}\rangle_{\partial\Omega^q} \\ 
    &-\langle \phi^h \{\{\omega^h \vb{u}^h\}\},\vb{n}\rangle_{\partial\Omega^q} - \langle \nabla \phi^h \cdot \vb{t}, \alpha[[\vb{F}]]\cdot \vb{n}^+\rangle_{\partial\Omega^q}\text{,}
\end{split}
\end{equation}
where we have used the definition of $\hat{G}$. Using the tangent boundary SBP property in section 3.5.5 this becomes
\begin{equation}
\begin{split}
    \langle \phi^h, \omega^h_t\rangle_{\Omega^q} -\langle \nabla_d \phi^h, \omega^h \vb{u}^h\rangle_{\Omega^q}  + \langle \phi^h \{\{\omega^h \vb{u}^h\}\},\vb{n}\rangle_{\partial\Omega^q} \\ + \langle \nabla \phi^h \cdot \vb{t}, \alpha[[\vb{F}]]\cdot \vb{n}^+\rangle_{\partial\Omega^q} = 0\text{,}\forall \phi^h \in S\text{,}
\end{split}
\end{equation}
demonstrating that the potential $G$ does not spuriously force the discrete {absolute} vorticity.

To show local {absolute} vorticity conservation we choose $\phi=1$ to yield
\begin{equation}
    \langle 1, \omega^h_t\rangle_{\Omega^q} = {-}\langle \{\{\omega^h \vb{u}^h\}\},\vb{n}\rangle_{\partial\Omega^q}\text{.}
\end{equation}
By construction $\{\{\omega^h \vb{u}^h\}\}$ is continuous across element boundaries and therefore {absolute} vorticity is locally conserved. Global conservation in a periodic domain follows from cancellation of $\{\{\omega^h \vb{u}^h\}\} \cdot \vb{n}$ across element boundaries.
 
\subsection{Energy conservation and stability}
Here we show that, for centred flux with $\alpha = \beta = 0$, the semi-discrete energy is conserved  and for  $\alpha, \beta > 0$ the semi-discrete energy is dissipated. We will obtain fully discrete energy stability by using a strong stability preserving time integrator \citep{shu1988efficient}.

Consider the element ${\Omega^q}$, we define the discrete energy $\mathcal{E}^q$ within the element as
\begin{equation}
\mathcal{E}^q = \tfrac{1}{2}\langle D^h\vb{u}^h, \vb{u}^h\rangle_{\Omega^q} + \tfrac{1}{2}\langle gD^h, D^h\rangle_{\Omega^q}\text{,}
\end{equation}
with the total energy given by the sum over all elements $\mathcal{E} = \sum_q \mathcal{E}^q$. Following the continuous form of the energy conservation analysis in section \ref{sec:Conserve_mass_vorticity_energy} the semi-discrete time evolution of $\mathcal{E}^q$ is given by
\begin{equation}
\begin{split}
	\mathcal{E}^q_t &= \langle D^h\vb{u}^h, \vb{u}^h_t\rangle_{\Omega^q} + \langle \tfrac{1}{2}\vb{u}^h\cdot\vb{u}^h, D^h_t\rangle_{\Omega^q} + \langle gD^h, D^h_t\rangle_{\Omega^q}\\
	 &= \langle \vb{F}^h, \vb{u}^h_t\rangle_{\Omega^q} + \langle G^h, D^h_t\rangle_{\Omega^q}.
  \end{split}
\end{equation}
Using  the weak forms of the mass and velocity  equations \eqref{eq:dg-velocity}--\eqref{eq:dg-depth} gives
\begin{equation}
\begin{split}
	\mathcal{E}^q_t &=  -\langle \vb{F}^h, \nabla_d G^h\rangle_{\Omega^q} + \langle \nabla_d G^h,  \vb{F}^h\rangle_{\Omega^q} - \langle \hat{G}, \vb{F}^h\cdot\vb{n} \rangle_{\partial\Omega^q} \\
 &- \langle G^h, (\hat{\vb{F}}^h - \vb{F}^h)\cdot\vb{n} \rangle_{\partial\Omega^q}\\
&= - \langle \hat{G}, \vb{F}^h\cdot\vb{n} \rangle_{\partial\Omega^q} - \langle G^h, (\hat{\vb{F}}^h - \vb{F}^h)\cdot\vb{n} \rangle_{\partial\Omega^q}
  \end{split}
\end{equation}
where we have used the fact that $\vb{F}^h\cdot \vb{k} \cross \vb{u}^h = 0$.
Note that $\vb{F}^h \cdot \vb{k} \cross \vb{u}^h = D^h\vb{u}^h \cdot \vb{k} \cross \vb{u}^h = 0$ only holds for the special case of diagonal mass matrices \citep{giraldo2020introduction}.

To prove energy stability, that is $\mathcal{E}_t = \sum_q \mathcal{E}^q_t \leq 0$, it suffices to show that the jump of the numerical energy flux across element boundaries is negative semi-definite, that is
\begin{equation}
    -\Big[\Big[ \vb{F}^h\hat{G} + G^h\big(
\hat{\vb{F}}^h - \vb{F}^h\big) \Big]\Big] \cdot\vb{n}^{+}\leq 0,
\end{equation}
for all points on the element interface. Using the expansion $[[ab]] = \{\{a\}\}[[b]] + [[a]]\{\{b\}\}$ the energy jump condition becomes
{\small
\begin{equation}
\begin{split}
 -\Big[\Big[ \vb{F}^h\hat{G} &+ G^h\big(
\hat{\vb{F}}^h - \vb{F}^h\big) \Big]\Big] \cdot\vb{n}^{+}= \\
&\big(\{\{G^h\}\} - \hat{G}\big)[[ \vb{F}^h]]\cdot\vb{n}^{+} + [[G^h]]\big(\{\{\vb{F}^h\}\} - \hat{\vb{F}}\big)\cdot\vb{n}^{+} \leq 0\text{.}
\end{split}
\label{eq:e-jump}
\end{equation}
}
For centred fluxes with $\alpha =0$, $\beta =0$, we have  $\hat{G} = \{\{G^h\}\}$ and $\hat{\vb{F}} = \{\{\vb{F}\}\}$ the jump  is $0$ and therefore the energy is conserved,
$$
-\Big[\Big[ \vb{F}^h\hat{G} + G^h\big(
\hat{\vb{F}}^h - \vb{F}^h\big) \Big]\Big] \cdot\vb{n}^{+} = 0 \implies \mathcal{E}_t =0.
$$
When $\alpha \ge 0$, $\beta \ge 0$, we have  
 $\hat{G} = \{\{G^h\}\} + \alpha[[\vb{F}^h]]\cdot\vb{n}^{+}$, $\hat{\vb{F}}\cdot\vb{n}^{+} = \{\{\vb{F}^h\}\}\cdot\vb{n}^{+} + \beta[[G^h]]$ with 
\begin{equation}
	-\Big[\Big[ \vb{F}^h\hat{G} + G^h\big(
\hat{\vb{F}}^h - \vb{F}^h\big) \Big]\Big] \cdot\vb{n}^{+}=-\alpha\big([[ \vb{F}^h]]\cdot\vb{n}^{+}\big)^2 - \beta[[G^h]]^2 \leq 0\text{.}
\end{equation}
The energy/entropy is dissipated $\mathcal{E}_t \le 0$.

\subsection{Discrete steady linear modes}

For the linearised RSWE with constant $f$ detailed in section \ref{sec:linear_geo_balance} we can preserve geostrophic balance by ensuring that the projection of any continuous stream function $\psi$ into the discrete space $S$ is in $S\cap H^1(\Omega)$, and that the resulting initial conditions satisfy $D^h \in S \cap H^1(\Omega)$ and $\vb{u}^h \in V \cap H(div)(\Omega)$ and remain in these subspsaces for the discrete linear system. For the linearised RSWE our method becomes
\begin{equation}
	\langle \vb{w}^h, \vb{u}^h_t\rangle_{\Omega^q} +  \langle \vb{w}^h, f\vb{k}\times\vb{u}^{h}+g\nabla D^h\rangle + \langle g(\hat{D} - D^h)\vb{w}, \vb{n}\rangle_{\partial\Omega^q} = 0\text{,}
\end{equation}
\begin{equation}
	\langle \phi^h, D^h_t\rangle_{\Omega^q} +  \langle \phi^h, H\nabla \cdot \vb{u}^h\rangle_{\Omega^q} + \langle H(\hat{\vb{u}} - \vb{u}^h), \vb{n}\rangle_{\partial\Omega^q}= 0\text{.}
\end{equation}
Analogous to the continuous case for any continuous stream function $\psi$ the discrete initial condition $D^h = -\tfrac{f}{g}\psi^h$, $\vb{u}^h = \nabla_d \cross \psi^h\vb{k}$ where $\psi^h = \Pi_{S}(\psi)$ is an exact discrete steady state. To prove this we show that for these particular initial conditions both the interior and boundary terms vanish. 

To see that the interior terms vanish we apply the div-curl cancellation property yielding 
\begin{equation}
    \nabla_d \cdot \vb{u}^h = \nabla_d \cdot \nabla_d \cross \psi^h\vb{k} = 0\text{,}
\end{equation}
and by construction 
\begin{equation}
\begin{split}
    \langle \vb{w}^h, f\vb{k}\cross\vb{u}^h + g\nabla_d D^h\rangle_{\Omega^q} &= \langle \vb{w}^h, f\vb{k}\cross\nabla_d \cross \psi^h \vb{k}- f\nabla_d \psi^h\rangle_{\Omega^q} \\ &= \langle \vb{w}^h, f\nabla_d \psi^h - f\nabla_d \psi^h\rangle_{\Omega^q} = 0\text{.}
    \end{split}
\end{equation}
To show that the boundary terms are $0$ we use the continuity property described in section 3.5.3 to show that $D^h$ and the normal components of $\vb{u}^h$ are continuous across element boundaries, and  therefore $\hat{D} - D^h = 0$ and $(\hat{\vb{u}} - \vb{u}^h)\cdot\vb{n} = 0$. To see that $D^h$ is continuous across element boundaries note that continuity of $\psi$ implies continuity of $\psi^h = \Pi_S(\psi)$, and therefore 
\begin{equation}
    D = -\dfrac{f}{g}\psi^h = -\dfrac{f}{g}\Pi_S(\psi) \in S\cap H^1(\Omega)\text{.}
\end{equation}
is continuous as well. Next we show that $\vb{u}$ is continuous in the normal direction. As $\psi^h$ is continuous it's derivative in the tangent direction $\big(\nabla^d \psi^h\big) \cdot \vb{t} $ is also continuous across element boundaries. The initial flow is given by
\begin{equation}
    \vb{u} = \nabla_d \cross \psi^h\vb{k} = \vb{k}\cross\nabla_d \psi^h\text{,}
\end{equation}
the flow in the normal direction is then
\begin{equation}
    \vb{u}\cdot \vb{n} = \vb{k}\cross\nabla_d \psi^h \cdot \vb{n} = \nabla_d \psi^h \cdot (\vb{n} \cross \vb{k}) = -\nabla_d \psi^h \cdot \vb{t}\text{,}
\end{equation}
and is therefore continuous across element boundaries. Therefore both the interior and element boundary terms are zero, and hence for any geostrophic mode there is a corresponding discrete solution which is exactly steady and will be preserved up to machine precision. We note that this result is dependent on the use of a collocated tensor product basis.

\section{Numerical results}

In this section we present the results of numerical experiments which verify our theoretical analysis. For all experiments we use an equi-angular cubed sphere mesh \citep{ranvcic1996global}, a {strong stability preserving} RK3 time integrator \citep{shu1988efficient}, and {polynomials of degree 3}. Except where otherwise specified, we use an adaptive time step with a fixed CFL of 0.8, {with the CFL defined as $\text{CFL} := \frac{c(2p + 1)\Delta t}{\Delta x}$} where $p$ is the polynomial order, $c$ is the fastest wave speed, $\Delta x$ is the element spacing, and $\Delta t$ is the time step size.

\subsection{Discrete steady states}

To verify that our scheme supports discrete steady modes we use the stream function $\psi = 0.1 \cos(\lambda)\cos(\theta)$, where $\lambda$ and $\theta$ are latitude and longitude respectively, with $6\times 5\times5$ elements, a spatially constant Coriolis parameter $f=8$, $H=0.2$, and $g=8$. Following the analysis in section 3.10 we use the initial condition $\vb{u}^h = \nabla \cross \Pi_S(\psi)\vb{k}$ and $D^h = -\tfrac{f}{g}\Pi_S(\psi)$ and integrate forward in time. Figure 5 shows the time $L^2$ difference from the initial conditions, demonstrating that the initial state is preserved up to machine precision.

\begin{figure}[!hbtp]
\begin{center}
     \includegraphics[width=0.8\textwidth]{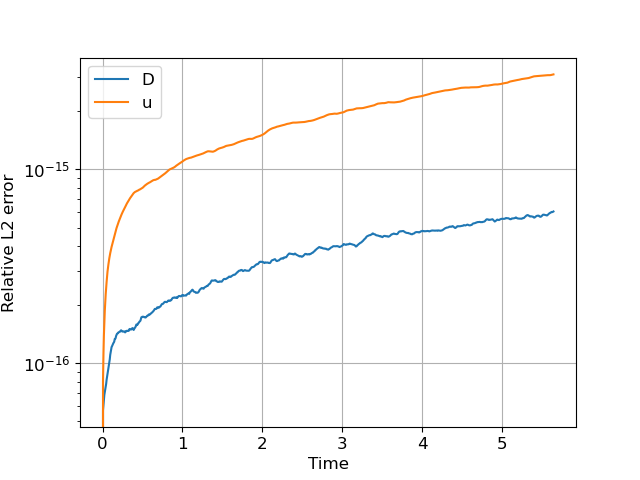}
     \caption{L2 difference from discrete initial conditions for the linear geostrophic balance test case.}
\end{center}
\end{figure}

\subsection{Williamson test case 2}

We apply both the energy conserving and energy dissipating methods to the steady state Williamson test case 2 with angle $\alpha=0$ \citep{williamson1992standard}. These initial conditions are
\begin{equation}
    u = u_0 \cos(\theta)\text{,}
\end{equation}
\begin{equation}
    D = D_0 - \dfrac{1}{g}\big(a f u_0 + \tfrac{1}{2}u_0^2\big)\sin^2(\theta)\text{,}
\end{equation}
where $\theta$ is latitude, $u$ is the zonal velocity, $a=6.37122\times 10^6$ m is the planetary radius, $g=9.80616 \text{ ms}^{-2}$, $u_0 = 38.61068 \text{ ms}^{-1}$ is the maximum velocity, $D_0 = \tfrac{1}{g}2.94\times 10^4$ m is the maximum depth, and $f = 2\times 7.292\times 10^{-5} \sin\theta \text{ s}^{-1}$.

We use a resolution of $6\times n^2$ elements for increasing $n=[3, 5, 10, 15, 30]$ {and calculate the error at day 5}. Figure 1 shows that both methods are converging monotonically at greater than 3rd order validating the consistency of our scheme. The energy conservative method is converging at approximately 3.4 order, while the dissipative method is converging at 3.8 order. This difference in convergence is most likely due to the energy dissipating method damping spurious high frequency waves. {We note that our methods do not attain the optimal 4th order convergence rate reported for the CG-SEM method in \citep{taylor2010compatible}, however we conjecture that this optimal 4th order convergence could be obtained by using similar entropy dissipative upwind fluxes to those in \citep{wintermeyer2017entropy, waruszewski2022entropy, winters2017uniquely}}.

\begin{figure}[!hbtp]
\begin{center}
	\includegraphics[width=0.7\textwidth]{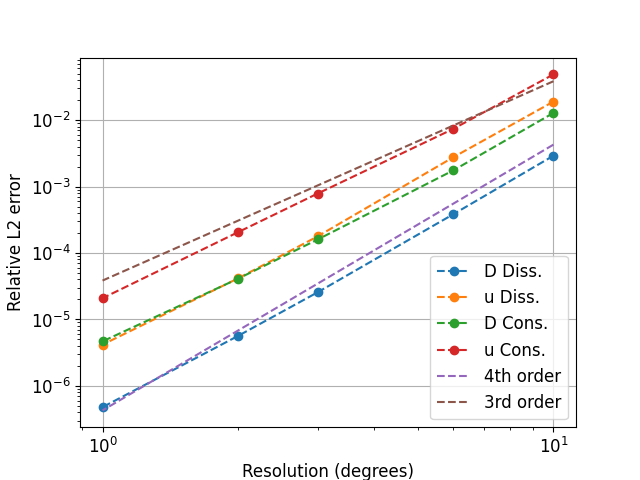}
 \caption{$L^2$ error for the energy dissipating and energy conserving methods applied to Williamson test case 2 {at day 5}.}
\end{center}
\end{figure}

\subsection{Williamson test case 5 - flow over an isolated mountain}

\begin{figure}[!hbtp]
\begin{center}
	\includegraphics[width=0.8\textwidth]{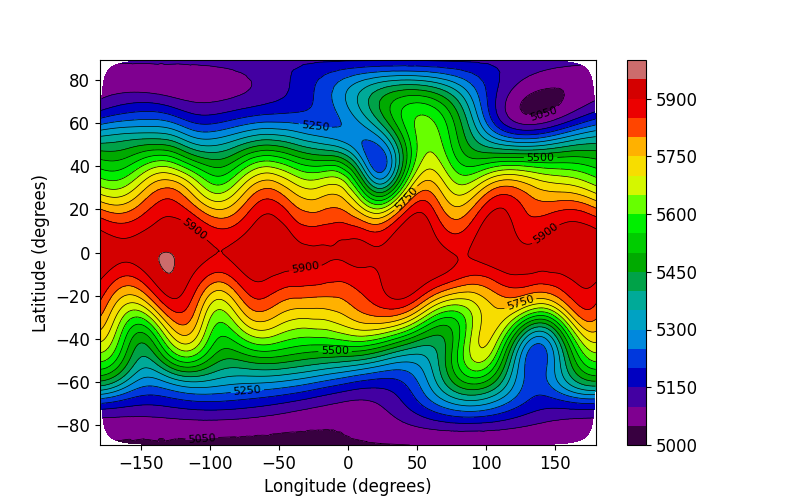}
 \caption{Test case 5 depth field at day 15 for the dissipative method at a resolution of $6\times32\times32$ elements. The contour interval is 50m.}
\end{center}
\label{fig:williamson-5}
\end{figure}

{
Test case 5 \citep{williamson1992standard} adds a single isolated mountain to test case 2 and uses $D_0=5960\text{ m}$, $u_0=20\text{ ms}^{-1}$, with all other parameters identical to test case 2. The bottom topography $b$ is given by
\begin{equation}
    b = \begin{cases}
        b_0(1 - \frac{r}{R}),&\text{for } r \leq R\\
        0, &\text otherwise
    \end{cases}
\end{equation}
where $r^2 = (\theta - \theta_c)^2 + \lambda^2$, $R=\frac{\pi}{9}$, $b_0 = 2000\text{ m}$, $\lambda$ is longitude, and $\theta_c=\frac{\pi}{6}$. We include the bottom topography in our scheme by adding the forcing term $-g\nabla_d b$ to the right hand side of \eqref{eq:dg-velocity}. Figure \ref{fig:williamson-5} displays the depth field at day 15 which closely agrees with the reference solution \citep{jakob1995spectral}.
}

\subsection{Williamson test case 6 (Rossby–Haurwitz wave)}

\begin{figure}[!hbtp]
\begin{center}
	\includegraphics[width=0.8\textwidth]{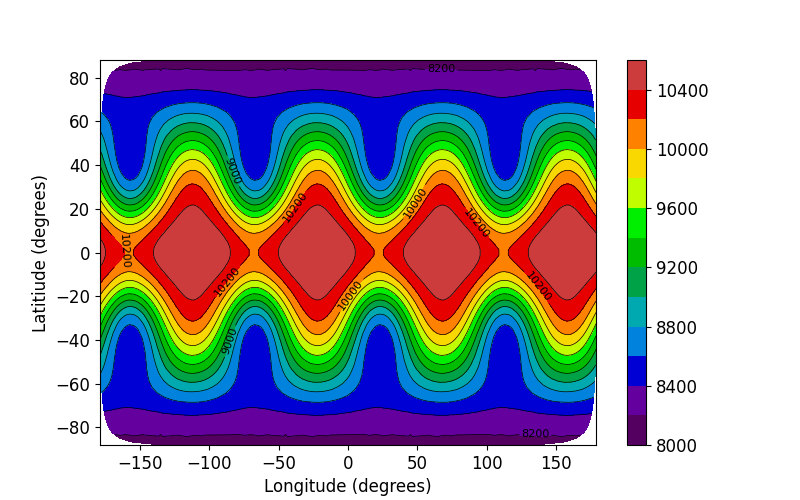}
 \caption{Test case 6 depth field at day 14 for the dissipative method at a resolution of $6\times16\times16$ elements. The contour interval is 200m.}
\end{center}
\label{fig:williamson-6}
\end{figure}

{
Test case 6 \citep{williamson1992standard} is a wavenumber 4 Rossby-Hauwitz wave which propogates west to east, without changing shape, with a group speed of $c_g = 7.7697\times10^7\text{ ms}^{-1}$ in the barotropic vorticity equations. This propogation is only approximate in the RSWE. We simulate this test case with the dissipative method at a resolution of $6\times16\times16$ elements for 14 days and find that the numerical group velocity magnitude is 94\% of $c_g$, in agreement with a previous numerical study \citep{lee2018mixed}. Figure \ref{fig:williamson-6} shows the depth field at day 14 which is consistent with other numerical results \citep{jakob1995spectral, taylor2010compatible, lee2018mixed}.
}

\subsection{Galewsky test case}

In order to validate the model in a more physically challenging turbulent regime we also simulate the Galewsky 
test case of a steady state zonal jet perturbed by a Gaussian hill in the depth field \citep{galewsky2004initial}. {The Gaussian hill induces gravity waves which excite a shear instability, causing the vorticity field to roll up into distinctive Kelvin-Helmholtz billows and the flow to become turbulent.} The initial conditions are 
\begin{equation}
    u = \begin{cases}
        \dfrac{u_0}{e_n} \exp\big((\theta-\theta_0)^{-1}(\theta-\theta_1)^{-1}\big)\text{, for } \theta_0 < \theta < \theta_1 \\ 
        0\text{, otherwise}
    \end{cases}
\end{equation}
\begin{equation}
    D = D_0 - \dfrac{a}{g}\int_{-\tfrac{\pi}{2}}^\theta{ u(\theta') \bigg[f + \dfrac{\tan(\theta')}{a}u(\theta')\bigg]d\theta'}\text{,}
\end{equation}
where $u$ is the zonal velocity, $\theta$ is latitude, $a=6.37122 \times 10^6$ m is the sphere's radius, $D_0 = 10^4$ m is the maximum depth, $u_0 = 80 \text{ ms}^{-1}$ is the maximum velocity, $e_n = \exp(-4 (\theta_1 - \theta_0)^{-2})$ is a normalisation constant, $\theta_0 = \tfrac{\pi}{7}$ and $\theta_1 = \tfrac{\pi}{2} - \theta_0$, $f=2 \times 7.292 \times 10^{-5} \sin\theta \text{ s}^{-1}$, and $g=9.80616 \text{ ms}^{-2}$.

\subsubsection{Discrete energy conservation}

To test the semi-discrete energy conservation of our energy conserving scheme we ran the Galewsky test case for 10 days at a low resolution of $6\times5\times5$ elements for decreasing time steps $\Delta t = [50s, 40s, 30s, 20s, 10s]$. {By day 10 the flow is turbulent, providing a challenging energy conservation test for our numerical scheme. Consistent with other energy conserving schemes \citep{taylor2010compatible, lee2018mixed, mcrae2014energy}} and as expected for a 3rd order time integrator, figure 2 shows that the energy conservation error decreases monotonically at 3rd order, verifying that the energy conserving scheme is semi-discretely energy conserving.

\begin{figure}[!hbtp]
\begin{center}
	\includegraphics[width=0.7\textwidth]{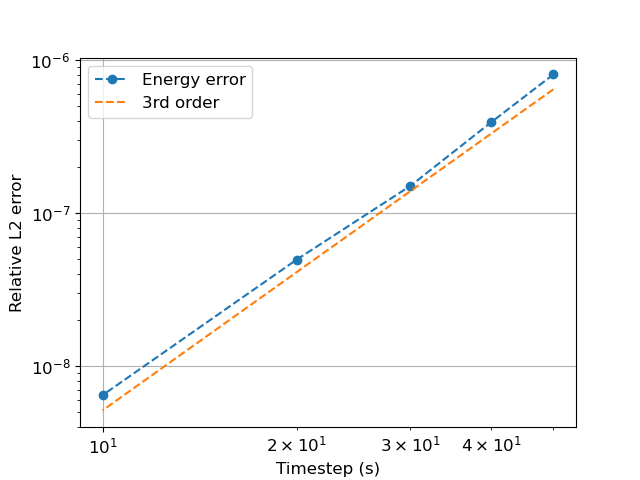}
 \caption{Energy conservation error for the centred flux scheme at low resolution for the Galewsky test case}
\end{center}
\end{figure}

\subsubsection{High resolution simulation}

To verify the accuracy and robustness of our method we ran the Galewsky test case for 20 days at a resolution of $6\times64\times64$ elements. Figure 3 shows that both the energy conserving and dissipating methods are energy stable, conserve mass and {absolute} vorticity   {with conservation errors consistent with other conservative schemes \citep{taylor2010compatible, lee2022comparison, mcrae2014energy}}, and appear to reasonably control potential enstrophy. It is noteworthy that our energy conserving method, which contains no dissipation of any kind, runs stably for the 20 day period. However, potential enstrophy is steadily increasing and the simulation may eventually become unstable. 

Figure 4 shows the relative vorticity at day 7 where the characteristic Kelvin-Helmholtz billows can be observed, {and their position, size, and shape closely match those of a previous study \citep{lee2022comparison}}. At day 7 some noise can be seen in the energy conserving solution, most likely due to spurious high frequency waves. {This noise is typical for numerical methods without dissipation \citep{galewsky2004initial, lee2022comparison}} but our entropy stable fluxes are able to damp these spurious numerical modes. 
To gain an indication as to whether the energy dissipating method is over damping we scale the final relative energy loss by the total energy in the Earth's atmosphere to obtain an adjusted average energy dissipation rate of $0.01Wm^{-2}$, an order of magnitude lower than the energy loss of the real atmosphere at the modelled length scales \citep{thuburn2008some}.

\begin{figure}[!hbtp]
\begin{center}
	\includegraphics[width=0.7\textwidth]{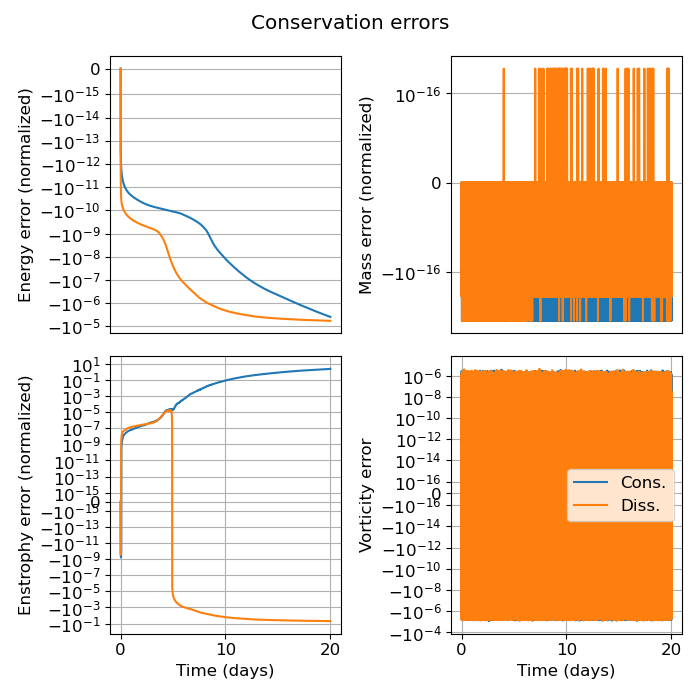}
 \caption{Conservation errors for the Galewsky test case at high resolution.}
\end{center}
\end{figure}

\begin{figure}[!hbtp]
\begin{center}
     \includegraphics[width=0.48\textwidth]{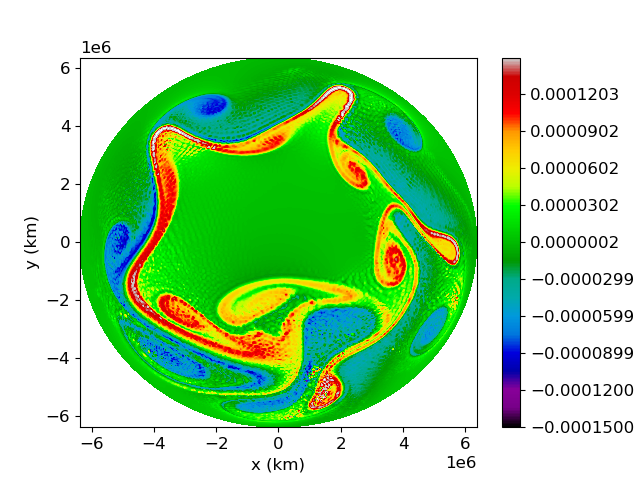}
     \includegraphics[width=0.48\textwidth]{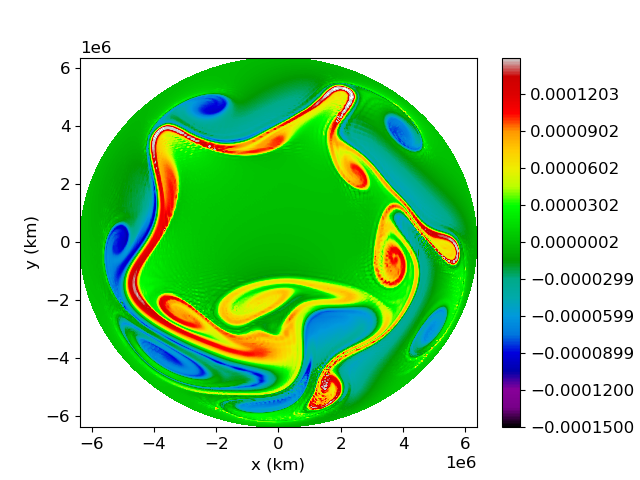} \\
     \caption{Relative vorticity for the Galewsky test case at day 7 for the energy conserving (left), and energy dissipating methods (right).}
     \label{fig:lin-geo}
\end{center}
\end{figure}

\section{Conclusion}

In this paper we presented a DG-SEM for the rotating shallow water equations in vector invariant form that possess many of the conservation properties and exact discrete linear geostrophic balance of recent mixed finite element methods. This enables our method to have the computational efficiency of SEM \citep{patera1984spectral} while having many of the mimetic properties of mixed finite element methods \citep{mcrae2014energy}.

We proved that on a general curvilinear grid our method locally conserves mass, {absolute} vorticity, and energy, does not spuriously project the potential onto {absolute} vorticity, and satisfies a discrete linear geostrophic balance. We then experimentally verified these properties and show that geostrophic turbulence can be well simulated with our entropy stable fluxes and no additional stabilisation.

While the issue of computational modes was not addressed, for the linear RSWE our method reduces to a DG method with Rusanov fluxes which is well known to have good dispersion properties \citep{rusanov1970difference}. Future work will focus on analysing the numerical modes of our method, investigating alternative numerical fluxes, and extending this method to the thermal shallow water equations.

\clearpage
\nocite{*}
\bibliographystyle{unsrtnat}
\bibliography{main}  

\end{document}